\documentclass[centertags,reqno]{amsart}
\usepackage{amsmath}
\usepackage{amsfonts}
\usepackage{amssymb}
\usepackage{graphicx}
\usepackage{amscd}
\usepackage{endnotes}
\usepackage{hyperref}

\theoremstyle{plain}

\numberwithin{equation}{section}

\begin{document}

\begin{center}
{\LARGE A multi-step approximant for fixed point problem and convex
optimization problem in Hadamard spaces\bigskip }

Muhammad Aqeel Ahmad Khan$^{\ast }$ and Hafiza Arham Maqbool

Department of Mathematics, COMSATS Institute of Information Technology,
Lahore 54000, Pakistan

February 27, 2018

\let\thefootnote\relax\footnotetext{%
* Corresponding author
\par
E-mail addresses: (MAA Khan)itsakb@hotmail.com,
maqeelkhan@ciitlahore.edu.pk, (A Maqbool) arhammaqbool101@gmail.com}
\end{center}

\textbf{Abstract}: The purpose of this paper is to propose and analyze a
multi-step iterative algorithm to solve a convex optimization problem and a
fixed point problem posed on a Hadamard space. The convergence properties of
the proposed algorithm are analyzed by employing suitable conditions on the
control sequences of parameters and the structural properties of the under
lying space. We aim to establish strong and $\triangle $-convergence results
of the proposed iterative algorithm and compute an optimal solution for a
minimizer of proper convex lower semicontinuous function and a common fixed
point of a finite family of total asymptotically nonexpansive mappings in
Hadamard spaces. Our results can be viewed as an extension and
generalization of various corresponding results established in the current
literature.\newline
\noindent \textbf{\noindent Keywords and Phrases}: Convex optimization,
Lower Semicontinuity, Proximal point algorithm, Total asymptotically
nonexpansive mapping, Common fixed point, Asymptotic center.\newline
\noindent \textbf{2010 Mathematics Subject Classification:} 47H09, 47H10,
65K10, 65K15.

\section{Introduction}

The theory of nonlinear analysis is mainly divided into three major areas,
namely convex analysis, monotone operator theory and fixed point theory of
nonlinear mappings. These theories have been largely developed in the
abstract setting of spaces having linear structures such as Euclidean,
Hilbert and Banach spaces. The theory of optimization, in particular, convex
optimization is prominent in the theory of convex analysis which studies the
properties of minimizers and maximizers of the under consideration
functions. The analysis of such properties rely on various mathematical
tools, topological notions and geometric ideas. Convex optimization not only
provides a theoretical setting for the existence and uniqueness of a
solution to a given optimization problem but also provides efficient
iterative algorithms to construct the optimal solution for such an
optimization problem. As a consequence, convex optimization solves a variety
of problems arising in disciplines such as mathematical economics,
approximation theory, game theory, optimal transport theory, probability and
statistics, information theory, signal and image processing and partial
differential equations, see, for example \cite%
{Adler,Combettes-2010,Combettes-2005,Smith-1994,Udriste} and the references
cited therein.

One of the major problems in optimization theory is to find a minimizer of a
convex function. The class of proximal point algorithms (PPA) contributes
significantly to the theory of convex optimization as to compute a minimizer
of a convex lower semicontinuous (lsc) function. In 1970, Martinet \cite%
{Martinet} proposed and analyzed the initial draft of PPA as a sequence of
successive approximation of resolvents. In 1976, Rockafellar \cite%
{Rockafellar} generally established, by the PPA, the convergence
characteristics to a zero of a maximal monotone operator in Hilbert spaces.
Brezis and Lions \cite{Brezis-Lions} improved the Rockafellar's algorithm
under a weaker condition on the parameters. The result established in \cite%
{Rockafellar} develops an interesting interplay between convex analysis,
monotone operator theory and fixed point theory of nonlinear mappings. As a
consequence, the PPA becomes an efficient tool for solving optimization
problems, fixed point problems, variational inequality problems and zeros of
maximal monotone operators. On the other hand, Rockafellar \cite{Rockafellar}
posed an open question regarding the strong convergence characteristics of
the PPA. The answer to the open question was settled in negative with a
counterexample given by G\"{u}ler \cite{Guler}. In order to establish strong
convergence of the PPA, one has to impose additional assumptions on the PPA,
see for example \cite{Boikanyo,Bruck-Reich,Kamimura,Xu-2006}. It is worth
mentioning that the counterexamples for strong convergence of the PPA are
still very rare and weak convergence is the best we can achieve without
additional assumptions.

Since most of the results in the theory of optimization involving PPA and
its various modifications are established within the spaces having linear
structure such as Euclidean space, Hilbert space and Banach space. It is
therefore natural to extend such beautiful and strong results from the
linear domain to the corresponding nonlinear domain. Another motivation for
this research direction is that various optimization problems, which are
non-convex in nature, become convex with the introduction of an adequate
metric defined on the under consideration spaces. Such metrics can also be
used to define new algorithms for optimization. Moreover, computation of
minimizers of the under consideration convex functions in such spaces plays
a pivotal role in the fields of nonlinear analysis and geometry \cite%
{Jost-1995,Jost-1997}. It is worth to mention that some efforts have been
made to generalize such results from the linear spaces to nonlinear spaces
having non-positive sectional curvature, see, for example, \cite%
{Bacak,Cholamjiak-2015,Cholamjiak,Ferreira,Li-Yao,Quiroz,Wang-Lopez} and the
references cited therein. This research area is still open either to
establish new convergence results for the class of PPA or to translate the
existing linear version of a result into the corresponding nonlinear version
in such spaces.

The outline of the paper is as follows: In Section 2, we first define the
conventions to be held throughout the paper and then define the consequent
notions, concepts and necessary results in the form of lemmas as required in
the sequel. Section 3 is devoted for the convergence analysis of the
proposed multi-step PPA to solve a convex optimization problem and a fixed
point problem posed on a Hadamard space.

\section{Preliminaries}

This section is devoted to recall some fundamental definitions, properties
and notations concerned with the fixed point problem and convex optimization
problem in Hadamard spaces. We also list some useful results in the form of
lemmas as required in the sequel. Throughout this paper, we write $%
x_{n}\rightarrow x~($resp. $x_{n}\rightharpoonup x)$ to indicate the strong
convergence (resp. the weak convergence) of a sequence $\{x_{n}\}_{n=1}^{%
\infty }$. The set of fixed points of a self-mapping $T$ on a nonempty
subset $C$ of a metric space $(X,d)$ is defined and denoted as: $F(T)=\{x\in
C:T(x)=x\}.$

Let $(X,d)$ be a metric space and $x,y\in X$ with $l=d\left( x,y\right) .$ A
geodesic from $x$ to $y$ in $X$ is a mapping $\theta :[0,l]\rightarrow X$
such that%
\begin{equation*}
\theta \left( 0\right) =x,\text{ }\theta \left( l\right) =y\text{ and }%
d\left( \theta (s),\theta (t)\right) =\left\vert s-t\right\vert ~\text{for
all}~s,t\in \lbrack 0,l].
\end{equation*}

The above characteristics shows that $\theta $ is an isometry and $x=\theta
(0)$ and $y=\theta (l)$ represent the end points of the geodesic segment.
The metric space $(X,d)$ is called a geodesic space if for every pair of
points $x,y\in X,$ there is a geodesic segment from $x$ to $y.$ Moreover, $%
(X,d)$ is uniquely geodesic if for all $x,y\in X$ there is exactly one
geodesic from $x$ to $y.$ A unique geodesic segment from $x$ to $y$ is
denoted as $[x,y].$ A geodesic triangle $\triangle \left(
x_{1},x_{2},x_{3}\right) $ in a geodesic metric space $(X,d)$ consists of
three points ${x_{1},x_{2},x_{3}}$ in $X$ (the vertices of $\triangle $) and
a geodesic segment between each pair of vertices (the edges of $\triangle $%
). A comparison triangle for the geodesic triangle $\triangle \left(
x_{1},x_{2},x_{3}\right) $ in $(X,d)$ is a triangle $\overline{\triangle }%
\left( x_{1},x_{2},x_{3}\right) :=\triangle \left( \overline{x}_{1},%
\overline{x}_{2},\overline{x}_{3}\right) $ in the Euclidean space $\mathbb{E}%
^{2}$ such that $d_{\mathbb{E}^{2}}\left( \bar{x_{i}},\bar{x_{j}}\right)
=d(x_{i},x_{j})$ for each $i,j\in \left\{ {1,2,3}\right\} .$

A geodesic space is said to be a $CAT(0)$ space if it is geodesically
connected and if every geodesic triangle in $(X,d)$ is at least as thin as
its comparison triangle in the Euclidean plane, that is $d(x,y)\leq d_{%
\mathbb{E}^{2}}(\bar{x},\bar{y}).$ A complete $CAT(0)$ space is then called
a Hadamard space. A nonempty subset $C$ of a $CAT(0)$ space is said to be
convex if $[x,y]\subset {C.}$ For a detailed discussion on this topic, we
refer the reader to consult \cite{Bridson-Haefliger,Bruhat-Tits}.

It is well known that a geodesic space is a $CAT(0)$ space if and only if%
\begin{equation*}
d^{2}((1-t)x\oplus t{y},z)\leq (1-t)d^{2}(x,z)+td^{2}(y,z)-t(1-t)d^{2}(x,y),
\end{equation*}%
\newline
for all ${x,y,z}\in {X}$ and ${t}\in \lbrack 0,1].$ In particular, if ${x,y}$
and ${z}$ are points in a $CAT(0)$ space and ${t}\in \lbrack 0,1],$ then%
\begin{equation*}
d((1-t)x\oplus t{y},z)\leq (1-t)d(x,z)+{t}d(y,z).
\end{equation*}

A self-mapping $T:C\rightarrow C$ is said to be total asymptotically
nonexpansive mapping \cite{Alber FPTA} if there exists non-negative real
sequences $\{k_{n}\}$ and $\{\varphi _{n}\}$ with $k_{n}\rightarrow 0$ and $%
\varphi _{n}\rightarrow 0$ as $n\rightarrow \infty $\ and a strictly
increasing continuous function $\xi :%
\mathbb{R}
^{+}\rightarrow
\mathbb{R}
^{+}$ with $\xi (0)=0$ such that
\begin{equation*}
d(T^{n}x,T^{n}y)\leq d(x,y)+k_{n}\xi \left( d(x,y)\right) +\varphi _{n}\text{
\ for all \ }x,y\in C,~n\geq 1.
\end{equation*}

The class of total asymptotically nonexpansive mappings is the most general
class of nonlinear mappings and contains properly various classes of
mappings associated with the class of asymptotically nonexpansive mappings.
These classes of nonlinear mappings have been studied extensively in the
literature \cite{FKK 2012,Khan JIA,Khan Fukhar FPTA,Khan-Fukhar-Amna} and
the references cited therein. It is worth mentioning that the results
established for total asymptotically nonexpansive mappings are applicable to
the mappings associated with the class of asymptotically nonexpansive
mappings and which are extensions of nonexpansive mappings.

It is well known that the concept of weak convergence in Hilbert spaces has
been generalized to $CAT(0)$ spaces as $\triangle $-convergence. Moreover,
many useful results from linear spaces involving weak convergence have
precise analogue version of $\triangle $-convergence in geodesic spaces. The
notion of asymptotic center of a sequence plays a key role to define the
concept of $\triangle $-convergence in such spaces.

Let $\{x_{n}\}$ be a bounded sequence in a $CAT(0)$ $X$. For $x\in X$,
define a continuous functional $r(\cdot ,\{x_{n}\}):X\rightarrow \lbrack
0,\infty )$ by:%
\begin{equation*}
r(x,\{x_{n}\})=\limsup_{n\rightarrow \infty }d(x,x_{n}).
\end{equation*}%
The asymptotic radius and asymptotic center of the bounded sequence $%
\{x_{n}\}$ with respect to a subset $C$ of $X$ is defined and denoted as:%
\begin{equation*}
r_{C}(\{x_{n}\})=\inf \{r(x,\{x_{n}\}):x\in C\},
\end{equation*}%
and%
\begin{equation*}
A_{C}(\{x_{n}\})=\{x\in C:r(x,\{x_{n}\})\leq r(y,\{x_{n}\})\text{ for all }%
y\in C\},
\end{equation*}%
respectively.

Recall that a sequence $\{x_{n}\}$ in $X$ is said to $\triangle $-converge
to $x\in X$ if $x$ is the unique asymptotic center of $\{u_{n}\}$ for every
subsequence $\{u_{n}\}$ of $\{x_{n}\}.$ In this case, we write $\triangle
-\lim_{n}x_{n}=x$ and call $x$ as the $\triangle $-limit of $\{x_{n}\}.$ A
mapping $T:C\rightarrow C$ is: (i) \textit{semi-compact} if every bounded
sequence $\{x_{n}\}\subset C$ satisfying $d(x_{n},Tx_{n})\rightarrow 0,$ has
a convergent subsequence; (ii) demiclosed at origin if for any sequence $%
\{x_{n}\}$ in $C$ with $x_{n}\rightharpoonup x$ and $\left\Vert
x_{n}-Tx_{n}\right\Vert \rightarrow 0,$ we have $x=Tx.$ Let $g$ \ be a
nondecreasing self-mapping on $[0,\infty )$ with $g(0)=0$ and $g(t)>0$\ for
all $t\in (0,\infty ).~$Let $\{{T}_{i}\}_{i=1}^{m}$ be a finite family of
total asymptotically nonexpansive mappings on $C$ with $\cap _{i=1}^{m}{F({%
T_{i}})\neq \emptyset }.$ Then the family of mappings is said to satisfy
Condition (I) on $C$ if:%
\begin{equation*}
\left\Vert x-Tx\right\Vert \geq f(d(x,F)),\text{ for all}\ x\in C,
\end{equation*}%
holds for at least one $T\in \{{T}_{i}\}_{i=1}^{m}.$

We now collect some basic concepts related to convex optimization in $CAT(0)$
spaces:

Let $C$ be a nonempty subset of a $CAT(0)$ space $X$, then a function $%
f:C\rightarrow (-\infty ,\infty ]$ is said be convex if for any geodesic $%
\theta :[a,b]\rightarrow C$ the function $f\circ {\theta }$ is convex. Some
important examples of convex function in $CAT(0)$ spaces can be found in
\cite{Bridson-Haefliger}. A function $f$ defined on $C$ is said to be lsc at
a point $x\in C$ if $f(x)\leq \liminf_{n\rightarrow \infty }f(x_{n}),$ for
each sequence $x_{n}\rightarrow x$. A function $f$ is said to be lsc on $C$
if it is lsc at any point in $C.$ A convex minimization problem associated
with a proper and convex function is to solve $x\in C$ such that%
\begin{equation*}
f(x)=\min_{y\in C}f(y).
\end{equation*}%
We denote by $\arg \min_{y\in C}f(y)$ by the set of a minimizer of a convex
function. For all $k>0$, define the \textit{Moreau-Yosida resolvent} of $f$
in a complete $CAT(0)$ space $X$ as follows:
\begin{equation*}
J_{k}(x)=\underset{y\in C}{\arg \min }[f(y)+\frac{1}{2k}{d^{2}(y,x)}],
\end{equation*}%
\newline
\newline
and put $J_{0}(x)=x$ for all $x\in {X.}$ This definition in metric spaces
with no linear structure first appeared in \cite{Guler}, see also \cite%
{Jost-1995}. The mapping $J_{k}$ is well defined for all $k\geq 0$ (see \cite%
{Guler,Jost-1995,Mayer}). For a proper, convex and lsc function, the set of
fixed points of the resolvent $J_{k}$ associated with $f$ coincides with the
set of minimizers of $f$\cite{Ariza-Leustean}. Moreover, the resolvent $%
J_{k} $ of $f$ is nonexpansive for all $k>0$ \cite{Jost-1995}. Some other
relevant characteristics of the resolvent $J_{k}$ of $f$ are incorporated in
the following couple of lemmas: \newline
\textbf{Lemma 2.1 }(Sub-differential Inequality) \cite{Ambrosio}. Let $(X,d)$
be a complete $CAT(0)$ space and $f:X\rightarrow (-\infty ,\infty ]$ be a
proper convex and lsc function. Then, for all $x,y\in {X}$ and {$k>0$}, we
have:%
\begin{equation*}
\frac{1}{2k}{d^{2}(J_{k}x,y)}-\frac{1}{2k}{d^{2}(x,y)}+\frac{1}{2k}{%
d^{2}(J_{k}x,x)}+f(J_{k}x)\leq f(y).
\end{equation*}%
\newline
\newline
\textbf{Lemma 2.2 }(The Resolvent Identity) \cite{Jost-1995,Mayer}. Let $%
(X,d)$ be a complete $CAT(0)$ space and $f:X\rightarrow (-\infty ,\infty ]$
be a proper convex and lsc function. Then, the following identity holds:%
\begin{equation*}
J_{k}x=J_{\mu }\left( \frac{k-\eta }{k}J_{k}x\oplus \frac{\eta }{k}{x}%
\right) .
\end{equation*}%
\newline
for all $x\in {X}$ and $k>\eta >0$.

We also require the following useful lemma for our main result.\newline
\textbf{Lemma 2.3 }\cite{Xu-2006}. Let $\left\{ a_{n}\right\} ,\left\{
b_{n}\right\} $\ and $\left\{ c_{n}\right\} $\ be sequences of non-negative
real numbers such that$~\sum_{n=1}^{\infty }b_{n}<\infty $ and $%
\sum_{n=1}^{\infty }c_{n}<\infty .$\ If~$a_{n+1}\ \leq (1+b_{n})a_{n}+c_{n},$%
\ $n\geq 1,~$then $\lim_{n\rightarrow \infty }a_{n}$\ exists\textit{.}
\newline
\textbf{Lemma 2.4 }\cite{KFK FPTA}. Let $(X,d,W)$\ be a uniformly convex
hyperbolic space with monotone modulus of uniform convexity $\eta .$~Let $%
x\in X$\ and $\{\alpha _{n}\}$\ be a sequence in $[a,b]$ for some $a,b\in
(0,1).$ If $\{x_{n}\}$ and $\{y_{n}\}$\ are sequences in $X$\ such that $%
\limsup\limits_{n\longrightarrow \infty }d(x_{n},x)\leq
c,~\limsup\limits_{n\longrightarrow \infty }d(y_{n},x)\leq c$\ and $%
\lim\limits_{n\rightarrow \infty }d(W(x_{n},y_{n},\alpha _{n}),x)=c$\ for
some $c\geq 0,$\ then $\lim\limits_{n\rightarrow \infty }d(x_{n},y_{n})=0.$%
\newline
\textbf{Lemma 2.5 (}\cite{KFK FPTA}\textbf{)}. \textit{Let }$K$\textit{\ be
a nonempty closed convex subset of a uniformly convex hyperbolic space and }$%
\{x_{n}\}$\textit{\ a bounded sequence in }$K$\textit{\ such that }$%
A_{K}(\{x_{n}\})=\{y\}$ and $r_{K}(\{x_{n}\})=\rho $\textit{. If }$\{y_{m}\}$%
\textit{\ is another sequence in }$K$\textit{\ such that }$%
\lim\limits_{m\rightarrow \infty }r(y_{m},\{x_{n}\})=\rho ,$ \textit{then }$%
\lim\limits_{m\rightarrow \infty }y_{m}=y.$\newline

\section{Main results}

We now prove a result in the form of lemma which plays a critical role to
establish strong and $\triangle $-convergence results of the proposed
iterative algorithm and compute an optimal solution for a minimizer of
proper convex lower semicontinuous function and a common fixed point of a
finite family of total asymptotically nonexpansive mappings in Hadamard
spaces.\newline
\textbf{Lemma 3.1.} Let $C$ be a nonempty closed convex subset of a Hadamard
space $X$. Let $f:X\rightarrow (-\infty ,\infty ]$ be a proper convex and
lsc function and let $\{{T}_{i}\}_{i=1}^{m}:C\longrightarrow C$ be a finite
family of uniformly continuous total asymptotically quasi nonexpansive
mappings with sequences $\{\lambda _{in}\}$ and $\{\mu _{in}\},$ $n\geq 1,\
i=1,2,\cdots ,m,$ such that\newline
(C1) $\sum\limits_{n=1}^{\infty }\lambda _{in}<\infty $ and$\
\sum\limits_{n=1}^{\infty }\mu _{in}<\infty ;$\newline
(C2) there exists constants $M_{i},\ M_{i}^{\ast }>0$ such that $\xi
_{i}\left( \theta _{i}\right) \leq M_{i}^{\ast }\theta _{i}$\ for all $%
\theta _{i}\geq M_{i}.$\newline
Let $\left\{ {x_{n}}\right\} $ be a sequence generated in the following
manner:%
\begin{equation}
\begin{pmatrix}
x_{1}\in C, \\
x_{n+1}=(1-\alpha _{n})y_{1n}\oplus \alpha _{n}T_{1}^{n}y_{1n}, \\
y_{1n}=(1-\alpha _{n})y_{2n}\oplus \alpha _{n}T_{2}^{n}y_{2n}, \\
\vdots \\
y_{in}=(1-\alpha _{n})y_{(i+n)n}\oplus \alpha _{n}T_{i+1}^{n}y_{(i+1)n}, \\
\vdots \\
y_{(m-1)n}=(1-\alpha _{n})z_{n}\oplus \alpha _{n}T_{m}^{n}z_{n}, \\
z_{n}=\underset{y\in C}{\text{ }\arg \min }[f(y)+\frac{1}{2k_{n}}{%
d^{2}(y,x_{n})}],\text{ }n\geq 1,%
\end{pmatrix}
\label{3.1}
\end{equation}%
\newline
where $\left\{ {\alpha _{n}}\right\} $ is a sequence in $[0,1]$ with ${0}<{a}%
<\alpha _{n}\leq {b}<{1}$ for all $n\geq 1$ and for some constant ${a,b}$ in
$(0,1).$ Assume that%
\begin{equation*}
\mathbb{F}=\left( \bigcap_{i=1}^{m}{F({T_{i}})}\right) \cap \underset{y\in C}%
{\text{ }\arg \min }f(y)\neq \emptyset ,
\end{equation*}%
then, we have the following:\newline
(i) $\lim_{n\rightarrow \infty }d(x_{n},p)$ exists for all $p\in \mathbb{F}$;%
\newline
(ii) $\lim_{n\rightarrow \infty }d(x_{n},z_{n})=0$;\newline
(iii) $\lim_{n\rightarrow \infty }d(T_{i}x_{n},x_{n})=0$, for each $%
i=1,2,\cdots ,m$.\newline
\textbf{Proof}: Let $p\in \mathbb{F}$, then $p={{T_{i}(p)}}$ for each $%
i=1,2,\cdots ,m$ and $f(p)\leq {f(y)}$ for all $y\in {C.}$ This implies that%
\begin{equation*}
f(p)+\frac{1}{2k_{n}}{d^{2}(p,p)}\leq f(y)+\frac{1}{2k_{n}}{d^{2}(y,p)},
\end{equation*}%
for each $y\in {C.}$ Hence $p=J_{k_{n}}(p)$ for each $n\geq 1.$\newline
\textbf{(i).} Now, we first show that $\lim_{n\rightarrow \infty }d(x_{n},p)$
exists. Since $z_{n}=J_{k_{n}}{x_{n}}$ and $J_{k_{n}}$ is nonexpansive,
therefore, we have%
\begin{equation}
d(z_{n},p)=d(J_{k_{n}}{x_{n}},J_{k_{n}}{p})\leq {d(x_{n},p)}.  \label{3.2}
\end{equation}%
It follows from (3.1) that%
\begin{eqnarray*}
d(y_{(m-1)n},p) &=&d((1-\alpha _{n})z_{n}\oplus \alpha _{n}T_{m}^{n}z_{n},p)
\\
&\leq &(1-\alpha _{n})d(z_{n},p)+{\alpha _{n}}d(T_{m}^{n}z_{n},p) \\
&\leq &(1-\alpha _{n})d(z_{n},p)+{\alpha _{n}}\{d(z_{n},p)+\lambda _{mn}\xi
_{m}(d(z_{n},p))+\mu _{mn}\}.
\end{eqnarray*}%
Since $\xi _{m}$ is an increasing function, therefore $\xi _{m}\left(
d\left( x_{n},p\right) \right) \leq \xi _{m}\left( M_{m}\right) $ for $%
d\left( x_{n},p\right) \leq M_{m}.$ Moreover $\xi _{m}\left( d\left(
x_{n},p\right) \right) \leq d\left( x_{n},p\right) M_{m}^{\ast }$ for $%
d\left( x_{n},p\right) \geq M_{m}$( by C2). In either case, we have%
\begin{equation*}
\xi _{m}\left( d\left( x_{n},p\right) \right) \leq \xi _{m}\left(
M_{m}\right) +d\left( x_{n},p\right) M_{m}^{\ast },
\end{equation*}%
where $M_{m},\ M_{m}^{\ast }>0.$ As a consequence, we get%
\begin{equation}
d(y_{(m-1)n},p)\leq (1+{\alpha _{n}}{\lambda _{mn}}{M_{mn}^{\ast }}%
)d(z_{n},p)+{\alpha _{n}}\lambda _{mn}\xi _{m}({M_{mn}})+{\alpha _{n}}\mu
_{mn}.  \label{3.3}
\end{equation}%
Let ${a_{1}}=\max \left\{ {\alpha _{n}},{\alpha _{n}}\xi _{m}\left( {M_{mn}}%
\right) ,{\alpha _{n}}{M_{mn}^{\ast }}\right\} >0,$ the estimate (3.3)
becomes%
\begin{equation}
d(y_{(m-1)n},p)\leq (1+{a_{1}}{\lambda _{mn}})d(x_{n},p)+{a_{1}}(\lambda
_{mn}+\mu _{mn}).  \label{3.4}
\end{equation}%
Again, reasoning in the aforementioned manner, it follows from (3.1) that%
\begin{eqnarray*}
d(y_{(m-2)n},p) &=&d\left( (1-\alpha _{n}){y_{(m-1)n}}\oplus \alpha
_{n}T_{m-1}^{n}{y_{(m-1)n}},p\right) \\
&\leq &(1-\alpha _{n})d({y_{(m-1)n}},p)+{\alpha _{n}}d(T_{m-1}^{n}{y_{(m-1)n}%
},p) \\
&\leq &(1-\alpha _{n})d({y_{(m-1)n}},p)+{\alpha _{n}}\{d({y_{(m-1)n}}%
,p)+\lambda _{(m-1)n}\xi _{m-1}(d({y_{(m-1)n}},p))+\mu _{(m-1)n}\} \\
&\leq &\left( 1+{\alpha _{n}}{\lambda _{(m-1)n}}{M_{(m-1)n}^{\ast }}\right)
d({y_{(m-1)n}},p)+{\alpha _{n}}\lambda _{(m-1)n}\xi _{m-1}({M_{(m-1)n}})+{%
\alpha _{n}}\mu _{(m-1)n}.
\end{eqnarray*}%
Utilizing (3.3) in the above estimate and simplifying the terms, we have%
\begin{eqnarray}
d(y_{(m-2)n},p) &\leq &\left( 1+{\alpha _{n}}\lambda _{mn}{M_{mn}^{\ast }}%
+\left( {\alpha _{n}}{M_{(m-1)n}^{\ast }}+{\alpha _{n}^{2}}\lambda _{mn}{%
M_{mn}^{\ast }}{M_{(m-1)n}^{\ast }}\right) \lambda _{(m-1)n}\right)
d(x_{n},p)  \notag \\
&&+{\alpha _{n}}\lambda _{mn}\xi _{m}({M_{mn}})+{\alpha _{n}^{2}}\lambda
_{mn}\lambda _{(m-1)n}\xi _{m}({M_{mn}}){M_{(m-1)n}^{\ast }}+{\alpha _{n}}%
\mu _{mn}  \notag \\
&&+{\alpha _{n}}\mu _{(m-1)n}+{\alpha _{n}^{2}}\lambda _{(m-1)n}{%
M_{(m-1)n}^{\ast }}\mu _{mn}+{\alpha _{n}}\lambda _{(m-1)n}\xi _{m-1}({%
M_{(m-1)n}}).  \label{3.5}
\end{eqnarray}%
Let ${c_{1},c_{2}>0}$ be such that $\lambda _{mn}\leq {c_{1}}$ and $\lambda
_{(m-1)n}\leq {c_{2}}$ for all $n\geq {1.}$ Then for ${\alpha _{n}}\leq {b,}$
the estimate (3.5) simplifies as%
\begin{eqnarray*}
d(y_{(m-2)n},p) &\leq &\left( 1+b\lambda _{mn}{M_{mn}^{\ast }}+\left( b{%
M_{(m-1)n}^{\ast }}+b{c_{1}}{M_{mn}^{\ast }}{M_{(m-1)n}^{\ast }}\right)
\lambda _{(m-1)n}\right) d(x_{n},p) \\
&&+b\lambda _{mn}\xi _{m}({M_{mn}})+b{c_{2}}\lambda _{mn}\xi _{m}({M_{mn}}){%
M_{(m-1)n}^{\ast }}+b\mu _{mn}+b\mu _{(m-1)n} \\
&&+b{c_{2}}{M_{(m-1)n}^{\ast }}\mu _{mn}+b\lambda _{(m-1)n}\xi _{m-1}({%
M_{(m-1)n}}).
\end{eqnarray*}%
Similarly, let ${a_{2}}=\max \{b,b{M_{mn}^{\ast }},b({{M}_{(m-1)n}^{\ast }}+{%
c_{1}}{{M}_{mn}^{\ast }}{{M}_{(m-1)n}^{\ast }}),b({{c}_{2}}\xi _{m}({M_{mn}})%
{M_{(m-1)n}^{\ast }}\newline
+\xi _{m}({M_{mn}})),b(1+{c_{2}}{M_{(m-1)n}^{\ast }}),b\xi _{m-1}({M_{(m-1)n}%
})\}>{0.}$ Then the above estimate becomes%
\begin{equation}
d(y_{(m-2)n},p)\leq \left( 1+{a}_{2}\sum\limits_{i=m-1}^{m}\lambda
_{in}\right) d(x_{n},p)+{a}_{2}\sum\limits_{i=m-1}^{m}(\lambda _{in}+\mu
_{in}).  \label{3.6}
\end{equation}%
Continuing in the similar fashion, for any ${m}\geq {1,}$ we have%
\begin{equation}
d({x_{n+1}},p)\leq \left( 1+{a}_{m}\sum\limits_{i=1}^{m}\lambda _{in}\right)
d(x_{n},p)+{a}_{m}\sum\limits_{i=1}^{m}(\lambda _{in}+\mu _{in}),
\label{3.7}
\end{equation}%
for some constant ${{a}_{m}}>{0.}$\newline
It now follows from (C1) and Lemma 2.3 that $\lim_{n\rightarrow \infty
}d(x_{n},p)$ exists for all $p\in \mathbb{F}$. This completes the proof of
part (i).\newline
\textbf{(ii).} In order to proceed for part (ii), we assume, without loss of
any generality, that%
\begin{equation}
\lim_{n\rightarrow \infty }d(x_{n},p)={r}\geq {0}.  \label{3.8}
\end{equation}%
Taking $\lim \sup $ on both sides of the estimate (3.2), we have%
\begin{equation}
\limsup\limits_{n\rightarrow \infty }d(z_{n},p)\leq {r}.  \label{3.9}
\end{equation}%
Consider the following variant of the estimate (3.7)%
\begin{equation*}
d({x_{n+1}},p)\leq \left( 1+{a}_{m}\sum\limits_{i=1}^{m}\lambda _{in}\right)
d(z_{n},p)+{a}_{m}\sum\limits_{i=1}^{m}(\lambda _{in}+\mu _{in}).
\end{equation*}%
Applying $\lim \inf $ on both sides of the above estimate, we get%
\begin{equation}
\liminf\limits_{n\rightarrow \infty }d(z_{n},p)\geq {r}.  \label{3.10}
\end{equation}%
The estimates (3.9) and (3.10) collectively imply that%
\begin{equation}
\lim\limits_{n\rightarrow \infty }d(z_{n},p)={r}.  \label{3.11}
\end{equation}%
Now, from Lemma 2.1, we have%
\begin{equation*}
\frac{1}{2k_{n}}[d^{2}(z_{n},p)-d^{2}(x_{n},p)+d^{2}(x_{n},z_{n})]\leq
f(p)-f(z_{n}).
\end{equation*}%
Since $f(p)\leq f(z_{n})$ for each $n\geq {1}$, it follows that%
\begin{equation*}
d^{2}(x_{n},z_{n})\leq d^{2}(x_{n},p)-d^{2}(z_{n},p).
\end{equation*}%
Utilizing (3.8) and (3.11), the above estimate implies that%
\begin{equation}
\lim_{n\rightarrow \infty }d(x_{n},z_{n})=0.  \label{3.12}
\end{equation}%
This completes the proof of part (ii).\newline
\textbf{(iii).} We now establish asymptotic regularity of the sequence $%
\left\{ {x_{n}}\right\} $ involving a finite family of uniformly continuous
total asymptotically quasi nonexpansive mappings.

Consider the following another variant of the estimate (3.7)
\begin{equation*}
d({x}_{n+1},p)\leq \left( 1+{a}_{m-1}\sum\limits_{i=1}^{m-1}\lambda
_{in}\right) d(y_{(m-1)n},p)+{a}_{m-1}\sum\limits_{i=1}^{m-1}(\lambda
_{in}+\mu _{in}).
\end{equation*}%
\newline
Taking $\lim \inf $ on both sides of the above estimate, we get%
\begin{equation}
\liminf\limits_{n\rightarrow \infty }d(y_{(m-1)n},p)\geq {r}.  \label{3.13}
\end{equation}%
Moreover, taking $\lim \sup $ on both sides of (3.4), we have%
\begin{equation}
\limsup\limits_{n\rightarrow \infty }d(y_{(m-1)n},p)\leq {r}.  \label{3.14}
\end{equation}%
\newline
Hence, by (3.13) and (3.14), we obtain%
\begin{equation}
\lim\limits_{n\rightarrow \infty }d(y_{(m-1)n},p)=\lim\limits_{n\rightarrow
\infty }d\left( (1-\alpha _{n})z_{n}\oplus \alpha
_{n}T_{m}^{n}z_{n},p\right) ={r}.  \label{3.15}
\end{equation}%
It follows from the definition of ${T_{m}}$ that $\lim \sup_{n\rightarrow
\infty }d(T_{m}^{n}{z_{n}},p)\leq {r.}$ Utilizing this fact together with
(3.9) and (3.15), it then follows from Lemma 2.3 that\newline
\begin{equation}
\lim\limits_{n\rightarrow \infty }d({z_{n}},T_{m}^{n}{z_{n}})={0}.
\label{3.16}
\end{equation}%
Now, observe the following variant of (3.7)%
\begin{equation*}
d({x}_{n+1},p)\leq \left( 1+{a}_{m-2}\sum\limits_{i=1}^{m-2}\lambda
_{in}\right) d(y_{(m-2)n},p)+{a}_{m-2}\sum\limits_{i=1}^{m-2}(\lambda
_{in}+\mu _{in}).
\end{equation*}%
Taking $\lim \inf $ on both sides of the above estimate, we get%
\begin{equation}
\liminf\limits_{n\rightarrow \infty }d(y_{(m-2)n},p)\geq {r}.  \label{3.17}
\end{equation}%
\newline
Also, taking $\lim \sup $ on both sides of the estimate (3.6), we have%
\begin{equation}
\limsup\limits_{n\rightarrow \infty }d(y_{(m-2)n},p)\leq {r}.  \label{3.18}
\end{equation}%
Hence, by (3.17) and (3.18), we obtain%
\begin{equation}
\lim\limits_{n\rightarrow \infty }d(y_{(m-2)n},p)=\lim\limits_{n\rightarrow
\infty }d((1-\alpha _{n})y_{(m-1),n}\oplus \alpha
_{n}T_{(m-1)}^{n}y_{(m-1)n},p)={r}.  \label{3.19}
\end{equation}%
Again, it follows from the definition of ${T_{m-1}}$ that $\lim
\sup_{n\rightarrow \infty }d(T_{m-1}^{n}y_{(m-1)n},p,p)\leq {r.}$ Utilizing
this fact together with (3.15) and (3.19), it then follows from Lemma 2.3
that%
\begin{equation*}
\lim\limits_{n\rightarrow \infty }d(y_{(m-1)n},T_{m-1}^{n}y_{(m-1)n})={0}.
\end{equation*}%
Continuing in the similar fashion, we have%
\begin{equation}
\lim\limits_{n\rightarrow \infty }d(y_{in},T_{in}^{n}y_{in})={0},\text{ for }%
i=1,2,\cdots ,m-1.  \label{3.20}
\end{equation}%
Note that $d({x}_{n+1},y_{1n})\leq b\cdot d(y_{1n},T_{1n}^{n}y_{1n}).$
Therefore, letting ${n}\rightarrow {\infty }$ and utilizing (3.20), we get%
\begin{equation}
\lim\limits_{n\rightarrow \infty }d({x}_{n+1},y_{1n})=0.  \label{3.21}
\end{equation}%
Moreover, $d(y_{in},y_{\left( i+1\right) n})\leq b\cdot d(y_{\left(
i+1\right) n},T_{\left( i+1\right) n}^{n}y_{\left( i+1\right) n}),$ for $%
i=1,2,\cdots ,m-2.$ Again, letting ${n}\rightarrow {\infty }$ and utilizing
(3.20), we get%
\begin{equation}
\lim\limits_{n\rightarrow \infty }d(y_{in},y_{\left( i+1\right) n})=0,~\text{%
for }i=1,2,\cdots ,m-2.  \label{3.22}
\end{equation}%
As a consequence of the estimates (3.21) and (3.22), we have%
\begin{equation}
\lim\limits_{n\rightarrow \infty }d({x}_{n},y_{in})=0\text{ for }%
i=1,2,\cdots ,m-1.  \label{3.23}
\end{equation}%
Now, observe that%
\begin{eqnarray*}
d(T_{m}^{n}{x_{n}},{x_{n}}) &\leq &d(T_{m}^{n}{x_{n}},T_{m}^{n}{x_{n}}{z_{n}}%
)+d(T_{m}^{n}{z_{n}},{z_{n}})+d({z_{n}},{x_{n}}) \\
&\leq &Ld({x_{n}},{z_{n}})+d(T_{m}^{n}{z_{n}},{z_{n}})+d({z_{n}},{x_{n}}).
\end{eqnarray*}%
\newline
Letting ${n}\rightarrow {\infty }$ in the above estimate and utilizing
(3.12) and (3.16), we have%
\begin{equation}
\lim\limits_{n\rightarrow \infty }d(T_{m}^{n}{x_{n}},{x_{n}})={0}.
\label{3.24}
\end{equation}%
Similarly%
\begin{eqnarray*}
d(T_{m-1}^{n}{x_{n}},{x_{n}}) &\leq &d(T_{m-1}^{n}{x_{n}},T_{m-1}^{n}{%
y_{\left( m-1\right) n}})+d(T_{m-1}^{n}{y_{\left( m-1\right) n}},{y_{\left(
m-1\right) n}})+d({y_{\left( m-1\right) n}},{x_{n}}) \\
&\leq &Ld({x_{n}},{y_{\left( m-1\right) n}})+d(T_{m-1}^{n}{y_{\left(
m-1\right) n}},{y_{\left( m-1\right) n}})+d({y_{\left( m-1\right) n}},{x_{n}}%
).
\end{eqnarray*}%
Letting ${n}\rightarrow {\infty }$ in the above estimate and utilizing
(3.20) and (3.23), we have%
\begin{equation}
\lim\limits_{n\rightarrow \infty }d(T_{m-1}^{n}{x_{n}},{x_{n}})={0}.
\label{3.25}
\end{equation}%
Continuing in the similar fashion, we get%
\begin{equation*}
\lim\limits_{n\rightarrow \infty }d(T_{i}^{n}{x_{n}},{x_{n}})={0}\text{ for }%
i=1,2,\cdots ,m.
\end{equation*}%
Now, utilizing the uniform continuity of $T_{i},$ the following estimate:%
\begin{equation*}
d\left( x_{n},T_{i}x_{n}\right) \leq d\left( x_{n},T_{i}^{n}x_{n}\right)
+d\left( T_{i}^{n}x_{n},T_{i}x_{n}\right)
\end{equation*}%
implies that%
\begin{equation}
\lim\limits_{n\rightarrow \infty }d\left( T_{i}x_{n},x_{n}\right) =0\text{
for }i=1,2,\cdots ,m.  \label{3.26}
\end{equation}%
This completes the proof. \newline
\textbf{Theorem 3.2.} Let $C$ be a nonempty closed convex subset of a
Hadamard space $X$. Let $f:X\rightarrow (-\infty ,\infty ]$ be a proper
convex and lsc function and let $\{{T}_{i}\}_{i=1}^{m}:C\longrightarrow C$
be a finite family of uniformly continuous total asymptotically quasi
nonexpansive mappings with sequences $\{\lambda _{in}\}$ and $\{\mu _{in}\},$
$n\geq 1,\ i=1,2,\cdots ,m,$ such that\newline
(C1) $\sum\limits_{n=1}^{\infty }\lambda _{in}<\infty $ and$\
\sum\limits_{n=1}^{\infty }\mu _{in}<\infty ;$\newline
(C2) there exists constants $M_{i},\ M_{i}^{\ast }>0$ such that $\xi
_{i}\left( \lambda _{i}\right) \leq M_{i}^{\ast }\lambda _{i}$\ for all $%
\lambda _{i}\geq M_{i}.$\newline
Let $\left\{ {x_{n}}\right\} $ be the sequence generated in (3.1) such that%
\begin{equation*}
\mathbb{F}=\left( \bigcap_{i=1}^{m}{F({T_{i}})}\right) \cap \underset{y\in C}%
{\text{ }\arg \min }f(y)\neq \emptyset .
\end{equation*}%
Then the sequence $\left\{ {x_{n}}\right\} \ \triangle $-converges to a
common element of $\mathbb{F}$. \newline
\textbf{Proof}: In fact, it follows from (3.12) and Lemma 2.2, that%
\begin{eqnarray*}
d(J_{k}{x_{n}},{x_{n}}) &\leq &d(J_{k}{x_{n}},{z_{n}})+d({z_{n}},{x_{n}}) \\
&\leq &d(J_{k}{x_{n}},J_{{k}_{{n}}}{x_{n}})+d({z_{n}},{x_{n}}) \\
&\leq &d\left( J_{k}{x_{n}},J_{k}\left( \frac{k_{n}-k}{k_{n}}%
J_{k_{n}}x_{n}\oplus \frac{k}{k_{n}}{x_{n}}\right) \right) +d({z_{n}},{x_{n}}%
) \\
&\leq &d\left( {x_{n}},\left( 1-\frac{k}{k_{n}}\right) J_{k_{n}}x_{n}\oplus
\frac{k}{k_{n}}{x_{n}}\right) +d({z_{n}},{x_{n}}) \\
&=&\left( 1-\frac{k}{k_{n}}\right) d({x_{n}},J_{k_{n}}x_{n})+d({z_{n}},{x_{n}%
}) \\
&\leq &\left( 1-\frac{k}{k_{n}}\right) d({x_{n}},{z_{n}})+d({z_{n}},{x_{n}})
\\
&\rightarrow &{0}\text{ as }n\rightarrow \infty .
\end{eqnarray*}%
Moreover, it follows from Lemma 3.1(i) that $\lim_{n\rightarrow \infty
}d(x_{n},p)$ exists for all $p\in \mathbb{F}$, hence $\{x_{n}\}$ is bounded
and has a unique asymptotic center, that is, $A_{C}(\{x_{n}\})=\{x\}.$ Let $%
\{u_{n}\}$ be any subsequence of $\{x_{n}\}$ such that $A_{C}(\{u_{n}\})=\{u%
\}$ and by Lemma 3.1(iii), we have $\lim_{n\rightarrow \infty
}d(T_{i}u_{n},u_{n})=0$ for $i=1,2,\cdots ,m.$ Next, we show that $u\in
\mathbb{F}.$ For each $i\in \{1,2,3,\cdots ,m\},$ we define a sequence $%
\{z_{n}\}$ in $K$ by $z_{j}=T_{i}^{j}u.$ In the presence of increasing
function $\xi _{i}$ and (C2), we calculate%
\begin{eqnarray*}
d(z_{j},u_{n}) &\leq
&d(T_{i}^{j}u,T_{i}^{j}u_{n})+d(T_{i}^{j}u_{n},T_{i}^{j-1}u_{n})+\cdots
+d(T_{i}u_{n},u_{n}) \\
&\leq &d\left( u,u_{n}\right) +\lambda _{in}\xi _{i}\left( d\left(
u,u_{n}\right) \right) +\mu
_{in}+\sum_{r=0}^{j-1}d(T_{i}^{r}u_{n},T_{i}^{r+1}u_{n}) \\
&\leq &\left( 1+\lambda _{in}M_{i}^{\ast }\right) d(u,u_{n})+\lambda
_{in}\xi _{i}\left( M_{i}\right) +\mu
_{in}+\sum_{r=0}^{j-1}d(T_{i}^{r}u_{n},T_{i}^{r+1}u_{n}).
\end{eqnarray*}%
Taking $\lim \sup $ on both sides of the above estimate and utilizing (3.9)
and the fact that each $T_{i}$ is uniformly continuous, we have%
\begin{equation*}
r(z_{j},\{u_{n}\})=\limsup_{n\rightarrow \infty }d(z_{j},u_{n})\leq
\limsup_{n\rightarrow \infty }d(u,u_{n})=r(u,\{u_{n}\}).
\end{equation*}%
This implies that $\left\vert r(z_{j},\{u_{n}\})-r(u,\{u_{n}\})\right\vert
\rightarrow 0$ as $j\rightarrow \infty .$ It follows from Lemma 2.5 that $%
\lim_{j\rightarrow \infty }T_{i}^{j}u=u.$ Again, utilizing the uniform
continuity of $T_{i},$ we have that $T_{i}(u)=T_{i}(\lim_{j\rightarrow
\infty }T_{i}^{j}u)=\lim_{j\rightarrow \infty }T_{i}^{j+1}u=u.$ From the
arbitrariness of $i,$ we conclude that $u$ is the common fixed point of $%
\{T_{i}\}_{i=1}^{m}.$ It remains to show that $x=u.$ In fact, uniqueness of
the asymptotic center implies that%
\begin{eqnarray*}
\limsup_{n\rightarrow \infty }d(u_{n},u) &<&\limsup_{n\rightarrow \infty
}d(u_{n},x) \\
&\leq &\limsup_{n\rightarrow \infty }d(x_{n},x) \\
&<&\limsup_{n\rightarrow \infty }d(x_{n},u) \\
&=&\limsup_{n\rightarrow \infty }d(u_{n},u).
\end{eqnarray*}%
This is a contradiction. Hence $x=u.$ This implies that $u$ is the unique
asymptotic center of $\{x_{n}\}$ for every subsequence $\{u_{n}\}$ of $%
\{x_{n}\}.$ This completes the proof.\newline
\textbf{Remark 3.3.} It is worth mentioning that the analogous weak
convergence result in Hilbert spaces for the sequence $\{x_{n}\}$ defined in
(3.1) can easily be obtained as a corollary of Theorem 3.2.\medskip

We now establish strong convergence characteristics of the sequence $%
\{x_{n}\}$ defined in (3.1) in a Hadamard space $X.$ We first give a
necessary and sufficient condition for the strong convergence of the
sequence (3.1).\medskip \newline
\textbf{Theorem 3.4} Let $C$ be a nonempty closed convex subset of a
Hadamard space $X$. Let $f:X\rightarrow (-\infty ,\infty ]$ be a proper
convex and lsc function and let $\{{T}_{i}\}_{i=1}^{m}:C\longrightarrow C$
be a finite family of uniformly continuous total asymptotically quasi
nonexpansive mappings with sequences $\{\lambda _{in}\}$ and $\{\mu _{in}\},$
$n\geq 1,\ i=1,2,\cdots ,m,$ such that\newline
(C1) $\sum\limits_{n=1}^{\infty }\lambda _{in}<\infty $ and$\
\sum\limits_{n=1}^{\infty }\mu _{in}<\infty ;$\newline
(C2) there exists constants $M_{i},\ M_{i}^{\ast }>0$ such that $\xi
_{i}\left( \lambda _{i}\right) \leq M_{i}^{\ast }\lambda _{i}$\ for all $%
\lambda _{i}\geq M_{i}.$\newline
Let $\left\{ {x_{n}}\right\} $ be the sequence generated in (3.1) such that%
\begin{equation*}
\mathbb{F}=\left( \bigcap_{i=1}^{m}{F({T_{i}})}\right) \cap \underset{y\in C}%
{\text{ }\arg \min }f(y)\neq \emptyset .
\end{equation*}%
Then the sequence $\left\{ {x_{n}}\right\} \ $converges strongly to a point
in $\mathbb{F}$\ if and only if $\lim \inf_{n\rightarrow \infty }dist(x_{n},%
\mathbb{F})=0,$ where $dist\left( x,\mathbb{F}\right) =\inf \left\{ d\left(
x,p\right) :p\in \mathbb{F}\right\} .$ \newline
\textbf{Proof}: The necessity of the conditions is obvious. Thus, we only
prove the sufficiency. It follows from Lemma 3.1(i) that the sequence $%
\left\{ d(x_{n},p)\right\} _{n=1}^{\infty }$\ converges. Moreover, $\lim
\inf_{n\rightarrow \infty }d(x_{n},\mathbb{F})=0$ implies that $%
\lim_{n\rightarrow \infty }d(x_{n},\mathbb{F})=0.$ This completes the
proof.\medskip \newline
\textbf{Theorem 3.5} Let $C$ be a nonempty closed convex subset of a
Hadamard space $X$. Let $f:X\rightarrow (-\infty ,\infty ]$ be a proper
convex and lsc function and let $\{{T}_{i}\}_{i=1}^{m}:C\longrightarrow C$
be a finite family of uniformly continuous total asymptotically quasi
nonexpansive mappings with sequences $\{\lambda _{in}\}$ and $\{\mu _{in}\},$
$n\geq 1,\ i=1,2,\cdots ,m,$ such that\newline
(C1) $\sum\limits_{n=1}^{\infty }\lambda _{in}<\infty $ and$\
\sum\limits_{n=1}^{\infty }\mu _{in}<\infty ;$\newline
(C2) there exists constants $M_{i},\ M_{i}^{\ast }>0$ such that $\xi
_{i}\left( \lambda _{i}\right) \leq M_{i}^{\ast }\lambda _{i}$\ for all $%
\lambda _{i}\geq M_{i}.$\newline
Let $\left\{ {x_{n}}\right\} $ be the sequence generated in (3.1) such that%
\begin{equation*}
\mathbb{F}=\left( \bigcap_{i=1}^{m}{F({T_{i}})}\right) \cap \underset{y\in C}%
{\text{ }\arg \min }f(y)\neq \emptyset .
\end{equation*}%
Assume that $\left\{ T_{i},J_{k}\right\} $\ satisfies Condition (I), then
the sequence $\left\{ {x_{n}}\right\} \ $converges strongly to a point in $%
\mathbb{F}$.\newline
\textbf{Proof}: It follows from Lemma 3.1(iii) that
\begin{equation*}
\lim_{n\rightarrow \infty }d(T_{i}x_{n},x_{n})=0\text{ for }i=1,2,\cdots ,m.
\end{equation*}%
Moreover, from Theorem 3.2, we have%
\begin{equation*}
\lim_{n\rightarrow \infty }d(J_{k}{x_{n}},{x_{n}})=0.
\end{equation*}%
Since $\left\{ T_{i},J_{k}\right\} $ satisfies Condition (I), so we have,
either%
\begin{equation*}
\lim_{n\rightarrow \infty }g(d(x_{n},\mathbb{F}))\leq \lim_{n\rightarrow
\infty }d(T_{i}x_{n},x_{n})=0,
\end{equation*}%
or%
\begin{equation*}
\lim_{n\rightarrow \infty }g(d(x_{n},\mathbb{F}))\leq \lim_{n\rightarrow
\infty }d(J_{k}{x_{n}},{x_{n}})=0,
\end{equation*}%
In both cases, it imply that $\lim_{n\rightarrow \infty }g(d(x_{n},\mathbb{F}%
))=0.$ Since $g$ is nondecreasing and $g(0)=0,$ we have $\lim_{n\rightarrow
\infty }d(x_{n},\mathbb{F})=0.$ Rest of the proof follows from Theorem 3.4
and is, therefore, omitted.\newline
\textbf{Remark 3.6.} It is remarked that the strong convergence
characteristics of the sequence $\{x_{n}\}$ defined in (3.1) in a Hadamard
space $X$ can also be established by utilizing the compactness condition of $%
C$ or $T(C).$ Moreover, one utilize the modified version of the
semi-compactness condition satisfied by a family of mappings. We further
remark that our results can be viewed as an extension and generalization of
various corresponding results established in the current literature. In
particular: (i). Theorems 3.2 generalizes the corresponding results in \cite[%
Theorem 3]{L-2017}, \cite[Theorem 3.2]{Cholamjiak} and \cite[Theorem 3.2]%
{P-2017}; (ii). Theorem 3.4 generalizes the corresponding results in \cite[%
Theorem 5]{L-2017} and \cite[Theorem 3.5]{Cholamjiak} and (iii). Theorem 3.5
generalizes the corresponding results in \cite[Theorem 3.6]{Cholamjiak} and
\cite[Theorem 3.4]{P-2017}.\newline
\textbf{Open Questions:} (i). Can we modify the sequence (3.1) involving
nonself-mapping in a Hadamard space $X.$ (ii). Can we modify the sequence
(3.1) in the form of a shrinking projection method for the strong
convergence results in a Hadamard space $X.$

\end{document}